\documentstyle{article}

\begin{document}
\baselineskip=18pt
\title{\hspace{11cm}{\small\bf IMPNWU-971109}\\
\vspace{2cm}
The nondynamical r-matrix structure for the elliptic $A_{n-1}$ Calogero-Moser
model}

\author{
 Bo-yu Hou$^{b}$ and Wen-li Yang$^{a,b}$ 
\thanks{e-mail :wlyang@nwu.edu.cn}
\thanks{Fax    :0086-029-8303511}
\\
\bigskip\\
$^{a}$ CCAST ( World Laboratory ), P.O.Box 8730 ,Beijing 100080, China\\
$^{b}$ Institute of Modern Physics, Northwest 
University, Xian 710069, China
\thanks{Mailing address}}

\maketitle

\begin{abstract}

In this paper, we construct a new Lax operator for the elliptic
$A_{n-1}$ Calogero-Moser model with general $n(2\leq n$) from the classical 
dynamical twisting, in which the corresponding r-matrix is  purely numeric 
(nondynamical one). The nondynamical r-matrix structure
of this Lax operator is obtained, which is elliptic $Z_n$-symmetric
r-matrix.

\end{abstract}
\section{Introduction}
A general description of classical completely integrable models of n
one-dimensional particles with two-body interactions $V(q_i-q_j)$ was
given in Ref.[24]. To each simple Lie algebra and choice of one of these
type interactions, one can associate a classically completely integrable
systems[5,13,23,24]. The most general form of the potential in such models is
so-called elliptic Calogero-Moser(CM) model with an elliptic interacting
potential. The various degenerations of this general system yield
rational CM model (type I in Ref.[5]), hyperbolic CM model
(type II in Ref.[5]) and trigonometric CM model (type III in Ref.[5]).
So, the study of the elliptic CM model is of great importance in the
completely integrable particle systems.

The Lax pair representation (Lax representation) of a completely integrable
system, which means that the equations in the problem can be formulated in
a Lax form,is the most effective way to show its integrability and construct
the complete set of integrals of the motion.The Lax representation and its
corresponding r-matrix structure for rational,hyperbolic and trigonometric
$A_{n-1}$ CM models were constructed by Avan et al[5]; The Lax
representation for the elliptic CM models was constructed by Krichever[22] and
the corresponding r-matrix structure was given by Sklyanin[30] and
Braden et al[10].There exists a specific feature that the r-matrix of the Lax
representation for these models turn out to be a dynamical one (i.e it
depends upon the dynamical variables) and satisfy a dynamical Yang-Baxter
equations[6,10,11,30]. Such structures also appear in the study 
Ruijsenaars-Schneider
models, which is known as the relativistic Calogero-Moser models and can
be related to the soliton systems of the affine Toda field theories
[11,12,27,28].
Moreover, such a dynamical  r-matrix structure is connected with the
hamiltonian reduction
of the cotangent bundle of Lie algebra for Calogero-Moser model and
the contangent bundle of Lie group for Ruijsenaars-schneider model
[2,16,31].
This greatly promotes the study of the classical ( and quantum)
dynamical r-matrices (and R-matrix).
A partial classification scheme has very recently been proposed
for the dynamical r-matrix obeying the particular version of the dynamical
Yang-Baxter equation[3,14,29].However, at this time one lacks a general
classifying scheme such as exists in the case of nondynamical classical
r-matrices thanks to Belavin and Drinfeld[9].

The other difficulties presented by the dynamical aspect of the r-matrix also
occur : {\bf I.} the fundamental Poisson algebra of Lax operator, which
structural constants are given by a dynamical r-matrix, is generally
speaking no longer closed (cf. the nondynamical one );
{\bf II.} To solve the quantization problem
and its geometrical interpretation is still an open problem[30]. So far, only
for one particular case---the spin generalization of the CM model--- a proper
algebraic setting (the Gervais-Neveu-Felder equation) was found[4] which
allows one to quantize the model.On the other hand, it is well-known that
the Lax representation for a completely integrable models is not unique.
The different
 Lax representation of a integrable system is conjugated each other(only
 for finite particles system, but for the field system it should be 
 transfered
 as gauge transformation from each other).But the
 corresponding r-matrix should be transfered as a $``$gauge" transformation
 (see Eq.(5)) which is the classical dynamically twisted relations[4]
 between
 r-matrix. So , to overcome the above difficulties which caused by the
 dynamical r-matrix may be that whether another Lax representation for the
 CM model which has a nondynamical r-matrix structure could be found.The plan
 of our work is to find such a $``$good" Lax representation for the elliptic
 $A_{n-1}$ CM model if possible.In our former work[18], we succeeded in
 constructing a new Lax operator (cf. Krichever's[22]) for the elliptic
 $A_{n-1}$ CM model
 with $n=2$ and showing that the corresponding r-matrix is a nondynamical one
 which is the classical eight-vertex r-matrix[20].In present paper,
 extending  our former work in Ref.[18], we construct a new Lax operator (cf.
 Krichever's) for the elliptic $A_{n-1}$ CM model with the general
 $n(2\leq N)$ which  be a $``$good" one in sense that the corresponding
 r-matrix is nondynamical one----the classical $Z_n$-symmetric r-matrix.

 The paper is organized as follows. In section 2, from the classical
 dynamical twisting , the condition that the `` good" Lax representation
 could exist is found.In section 3, after some review of quantum $Z_n$-
 symmetry Belavin model, we construct the classical $Z_n$-symmetry r-matrix.
After  some review of Sklyanin's work on elliptic CM model in section 4, we
construct the ``good"  Lax  representation for elliptic $A_{n-1}$ CM model
which enjoys in nondynamical r-matrix structure in section 5. 
Finally, we give
summary and discussions in section 6. Appendix contains some detailed
calculations.

\section{The dynamical twisting of classical r-matrix}
In this paper, we only deal with the completely integrable finite
particles systems.In this section we will review some general theories of
the completely integrable finite particles systems

A Lax pair (L,M) consists of two functions on the phase space of the system
with values in some Lie algebra {\bf $g$}, such that the evolution equations may be
written in the following form
\begin{eqnarray}
\frac{dL}{dt}=[L,M]
\end{eqnarray}
\noindent where $[,]$ denotes the bracket in the Lie algebra {\bf $g$}.The
interest in the existence of such a pair lies in the fact that it allows
for an easy construction of conserved quantities (integrals of motion)---
it follows that the adjoint-invariant quantities $trL^{n}$ are the integrals
of the motion.In order to implement Liouville theorem onto this set of
possible action variables we need them to be Poisson-commuting.As shown in
Ref.[7], for the commutativity of the integrals $trL^{n}$ of the Lax
operator it is neccessary and sufficient that the fundamental Poisson bracket
$\{L_{1}(u),L_{2}(v)\}$ could be represented in the commutator form
\begin{eqnarray}
\{L_1(u),L_2(v)\}=[r_{12}(u,v),L_1(u)]-[r_{21}(v,u),L_{2}(v)]
\end{eqnarray}
\noindent where we use the notation
\begin{eqnarray*}
L_1\equiv L\otimes 1\ \ \ \ ,\ \ \ \ L_2\equiv 1\otimes L\ \ \ ,\ \ \
r_{21}=Pr_{12}P
\end{eqnarray*}
\noindent and $P$ is the permutation operator such that
$Px\otimes y=y\otimes x$.

Generally speaking, r-matrix $r_{12}(u,v)$ does depend on  dynamical
variables.For some special case of $r_{12}(u,v)$ independent on dynamical
variables,the r-matrix is called as the nondynamical r-matrix which has  been
 well studied[15].In contrast  with the well-studied case of the nondynamical
 r-matrix , no general theory of the dynamical r-matrix exists at the moment,
 apart few concrete examples and observations. Still, the collection
 of examples is rather sparse, and any new example of dynamical r-matrix
 could possibly contribute to better understanding of their algebraic and
 geometric nature.

The Poisson bracket structure Eq.(2) obeys a Jacobi identity which implies
an algebraic constraint for the r-matrix. Since r-matrix maybe depend on
dynamical variables this constraint takes a complicated form
\begin{eqnarray}
[L_1,[r_{12},r_{13}]+[r_{12},r_{23}]+[r_{32},r_{13}]+\{L_2,r_{13}\}-
\{L_3,r_{12}\}]+cycl.perm=0
\end{eqnarray}
\noindent Relevant particular cases of this general identity, one can obtain
the classical Yang-Baxter equation for the nondynamical r-matrix and the
classical dynamical Yang-Baxter equation [3,7]for the dynamical one.
It should be remarked that such a classification is by no means unique, which
drastically depend on the Lax representation that one choose for a system. 
Namely,
there is no one-to-one correspondence between a given dynamical system and
a defined r-matrix, a same dynamical system may have several Lax
representations and several r-matrix .The different Lax
representation of a system is conjugated each other: if ($\stackrel{\sim}{L},
\stackrel{\sim}{M}$) is one of other Lax pair of the same dynamical system
conjugated with the old one $(L,M)$, it means that 
\begin{eqnarray*}
\frac{d\stackrel{\sim}{L}}{dt}=[\stackrel{\sim}{L},\stackrel{\sim}{M}]
\end{eqnarray*}
\begin{eqnarray}
\stackrel{\sim}{L}(u)=g(u)L(u)g^{-1}(u)\ \ ,\ \
\stackrel{\sim}{M}(u)=g(u)M(u)g^{-1}(u) -(\frac{d}{dt}g(u))g^{-1}(u)
\end{eqnarray}
\noindent where $g(u)\in G$ whose Lie algebra is {\bf $g$}. Then, we have

\noindent {\bf Proposition 1.} The Lax pair
($\stackrel{\sim}{L},\stackrel{\sim}{M})$ has the following r-matrix
structure
$$
\{\stackrel{\sim}{L}_1(u),\stackrel{\sim}{L}_2(v)\}=
[\stackrel{\sim}{r}_{12}(u,v),\stackrel{\sim}{L}_1(u)]
-[\stackrel{\sim}{r}_{21}(v,u),\stackrel{\sim}{L}_{2}(v)]\eqno(4a)
$$
\noindent where
\begin{eqnarray}
& &\stackrel{\sim}{r}_{12}(u,v)=g_1(u)g_2(v)r_{12}(u,v)g^{-1}_1(u)g^{-1}_2(v)
+g_2(v)\{g_1(u),L_2(v)\}g^{-1}_1(u)g^{-1}_2(v)\nonumber\\
& &\ \ \ \ \ \ \ \ \ \ \ \ \ \ \ 
+\frac{1}{2}[\{g_1(u),g_2(v)\}g^{-1}_1(u)g^{-1}_2(v)\ \ ,
\ \ g_2(v)L_2(v)g^{-1}_2(v)]
\end{eqnarray}
\vspace{1cm}

\noindent {\it \bf Proof:} The proof is direct substituting 
Eq.(4) and Eq.(4a) into the
fundamental Poisson bracket and use the following identity
\begin{eqnarray*}
\left[ [s_{12},L_1],L_2\right] =\left[ [s_{12},L_2],L_1\right]
\end{eqnarray*}
\noindent where $s_{12}$ is any matrix on {\bf $g\otimes g$}.
\hspace{3cm} ${\bf\large\Box}$

It can be seen that: {\bf I.}The Lax operator $L$ is transfered
as a similarity transformation from the different Lax representation( only
for finite particles system);
{\bf II.}The corresponding $M$ is undergone the usual gauge transformation;
{\bf III.} The r-matrix is transfered as some generalized gauge
transformation, which can be considered as the classical version of the
dynamically twisted relation between the quantum R-matrix[4].Therefore,
it is of great value to find a $``$good" Lax
representation for a system if it exists, in which the corresponding 
r-matrix is nondynamical one and the well-studied
theories[4,15] can be directly applied in the system---such as the dressing
transformation , quantization.....

\noindent {\bf Corollary to proposition 1. } For given Lax pair $(L,M)$ and
the corresponding r-matrix, if there exist such a $g$ that
\begin{eqnarray}
& &h_{12}=\{g_1(u)g_2(v)r_{12}(u,v)g^{-1}_1(u)g^{-1}_2(v)
+g_2(v)\{g_1(u),L_2(v)\}g^{-1}_1(u)g^{-1}_2(v)\nonumber\\
& &\ \ \ \ \ \ \ \ \ \ \ 
+\frac{1}{2}[\{g_1(u),g_2(v)\}g^{-1}_1(u)g^{-1}_2(v)\ \ ,
\ \ g_2(v)L_2(v)g^{-1}_2(v)]\}\nonumber\\
& & {\rm and\ \ \  }\partial_{q_i}h_{12}=\partial_{p_j}h_{12}=0
\end{eqnarray}
\noindent the nondynamical Lax representation of the system exist.

\vspace{1cm}

\noindent From the straightforward calculation , we also have

\noindent {\bf Proposition 2.} The twisted Lax pair
($\stackrel{\sim}{L},\stackrel{\sim}{M})$ and the corresponding r-matrix
$\stackrel{\sim}{r}_{12}$ satisfies
\begin{eqnarray*}
[\stackrel{\sim}{L}_1,[\stackrel{\sim}{r}_{12},\stackrel{\sim}{r}_{13}]
+[\stackrel{\sim}{r}_{12},\stackrel{\sim}{r}_{23}]+
[\stackrel{\sim}{r}_{32},\stackrel{\sim}{r}_{13}]
+\{\stackrel{\sim}{L}_2,\stackrel{\sim}{r}_{13}\}-
\{\stackrel{\sim}{L}_3,\stackrel{\sim}{r}_{12}\}]+cycl.perm=0
\end{eqnarray*}
\vspace{1cm}

\noindent The main purpose of this paper is  to find a 
``good" Lax representation for the elliptic $A_{n-1}$ CM model.

\section{The elliptic function and elliptic $Z_n$-symmetric R-matrix and
r-matrix}
We first briefly review the elliptic $Z_n$-symmetric quantum R-matrix which
is related to $Z_n$-symmetric Belavin model[8,19,21,25].
For $n\in Z_{+},2\leq n$, we
define $n\times n$ matrices $h,g,I_{\alpha}$ as
\begin{eqnarray*}
h_{ij}=\delta_{i+1,j{\rm mod n}}\ \ ,\ \ g_{ij}=\omega^{i}\delta_{i,j}
\ \ ,\ \ I_{\alpha_1,\alpha_2}\equiv I_{\alpha}=g^{\alpha_2}h^{\alpha_1}
\end{eqnarray*}
\noindent where $\alpha_1,\alpha_2 \in Z_n$ and
$\omega=exp(2\pi\frac{\sqrt{-1}}{n})$. We also define
some elliptic functions
\begin{eqnarray}
& &\theta^{(j)}(u)=
\theta\left[\begin{array}{c}\frac{1}{2}-\frac{j}{n}\\ 
\frac{1}{2}\end{array}\right](u,n\tau)\ \ ,\ \
\sigma(u)=\theta\left[\begin{array}{c}\frac{1}{2}\\
\frac{1}{2}\end{array}\right](u,\tau)\\
& &\theta\left[\begin{array}{c}a\\ b\end{array}\right](u,\tau)
=\sum_{m=-\infty}^{\infty}exp\{\sqrt{-1}\pi[(m+a)^{2}\tau + 2(m+a)(z+b)]\}
\nonumber\\
& &\theta'^{(j)}(u)=\partial_{u}\{\theta^{(j)}(u)\}\ \ ,\ \
\sigma'(u)=\partial_{u}\{\sigma(u)\}\ \ ,\ \
\xi(u)=\partial_{u}\{ln\sigma(u)\}\\
& &E(u,v)=\frac{\sigma(u+v)}{\sigma(u)\sigma(v)}
\end{eqnarray}
\noindent where $\tau$ is a complex number with $Im(\tau)>0$ .
Then we define the $Z_n$-symmetric Belavin's R-matrix[19]
\begin{eqnarray}
R^{lk}_{ij}(v)=\left\{ \begin{array}{ll}
\frac{\theta^{'(0)}(0)\sigma(v)\sigma(w)}
{\sigma'(0)\theta^{(0)}(v)\sigma(v+w)}
\frac{\theta^{(0)}(v)\theta^{(i-j)}(v+w)}
{\theta^{(i-l)}(w)\theta^{(l-j)}(v)}&{\rm if }\ \ i+j=l+k\ \ {\rm mod\ \ n} \\
0&{\rm\ \ \ otherwise}
\end{array}
\right.
\end{eqnarray}
\noindent where $w$ is a complex number which is called as the crossing
parameter of the R-matrix.We should remark that our  R-matrix coincide
with the usual one [19,21] up to a scalar factor $
\frac{\theta^{'(0)}(0)\sigma(v)}{\sigma'(0)\theta^{(0)}(v)}
\prod_{j=1}^{n-1}\frac{\theta^{(j)}(v)}{\theta^{(j)}(0)}$
, which is to make Eq.(15) satisfied.The R-matrix satisfy quantum Yang-Baxter
equation (QYBE)
\begin{eqnarray}
R_{12}(v_1-v_2)R_{13}(v_1-v_3)R_{23}(v_2-v_3)=R_{23}(v_2-v_3)R_{13}(v_1-v_3)
R_{12}(v_1-v_2)
\end{eqnarray}
\noindent Moreover, the R-matrix enjoys in following $Z_n\otimes Z_n$
symmetric properties
\begin{eqnarray}
R_{12}(v)=(a\otimes a)R_{12}(v)(a\otimes a)^{-1}\ \ \ ,\ \ {\rm for}\ \ a=g,h
\end{eqnarray}

Introduce an $n\otimes n$ matrix $\widehat{T}(v)$, with the matrix element
$\widehat{T}(v)^{j}_{i}$ being operators,which satisfy the equation( also
called QYBE)
\begin{eqnarray}
R_{12}(v_1-v_2)\widehat{T}_1(v_1)\widehat{T}_2(v_2)
=\widehat{T}_2(v_2)\widehat{T}_1(v_1)R_{12}(v_1-v_2)
\end{eqnarray}
We next turn to the factorized difference representation for the operator
$\widehat{T}(v)$[17,19,25]

Set an $n\otimes n$ matrix $A(u;q)$
\begin{eqnarray}
& &A(u;q)^{i}_{j}\equiv A(u;q_1,q_2,...,q_n)^{i}_{j}=
\theta^{(i)}(u+nq_j-\sum_{k=1}^{n}q_k+\frac{n-1}{2})
\end{eqnarray}
\noindent here $A(u,q)^{i}_{j}$ correspond to the interwiner function
$\varphi^{(i)}_j$ between the $Z_n$-symmetric Belavin R-matrix and the
$A^{(1)}_{n-1}$ face model[21] in Ref.[19].
Construct the operator $\widehat{T}(u)$
$$
\widehat{T}(u)^{i}_{j}=A(u+sw;q)^{i}_{k}A^{-1}(u;q)^{k}_{j}D_{k}\eqno(14a)
$$
\noindent where $s$ is a complex number which associate with the
representation  of Sklyanin algebra[19] and will be related to the
coupling constant of elliptic $A_{n-1}$ CM model Eq.(36), and 
 $D_k$ is a difference operator such as
\begin{eqnarray*}
D_kf(q)\equiv D_kf(q_1,q_2,...,q_n)=f(q_1,...,q_{k-1},q_k-w,q_{k+1},...,q_n)
\end{eqnarray*}
Then following the results in Ref.[16], we have

\noindent {\bf  Theorem 1.} ([16],[29],[30]) The L-operator $\widehat{T}(u)$
defined in Eq.(14a) satsifes the QYBE Eq.(13).
\vspace{1cm}

\noindent We can define  a corresponding $Z_n$-symmetric (classical)
r-matrix which has the following relation with the R-matrix
\begin{eqnarray}
& &R_{12}(v)|_{w=0}=1\otimes 1\nonumber\\
& & R_{12}(v)=1\otimes 1 +wr_{12}(v)
+0(w^2)\ \ ,\ \ {\rm when\ \ the \ \ crossing\ \ parameter}\ \ w
\longrightarrow 0
\end{eqnarray}
\noindent Then we have

\noindent{\bf Proposition 3.} The corresponding elliptic $Z_n$-symmetric
r-matrix is 
\begin{eqnarray}
r^{lk}_{ij}(v)=\left\{\begin{array}{ll}
(1-\delta^{l}_{i})\frac{\theta^{'(0)}(0)\theta^{(i-j)}(v)}
{\theta^{(l-j)}(v)\theta^{(i-l)}(0)}
+\delta^{l}_{i}\delta^{k}_{j}(
\frac{\theta^{'(i-j)}(v)}{\theta^{(i-j)}(v)}-\frac{\sigma'(v)}{\sigma(v)})
&{\rm if }\ \ i+j=l+k\ \ {\rm mod}\ \ n\\
0&{\rm otherwise}
\end{array}
\right.
\end{eqnarray}
\noindent and satisfy the nondynamical (classical) Yang-Baxter equation and 
antisymmetric properties
\begin{eqnarray}
& &[r_{12}(v_1-v_2),r_{13}(v_1-v_3)]+[r_{12}(v_1-v_2),r_{23}(v_2-v_3)]
+[r_{13}(v_1-v_3),r_{23}(v_2-v_3)]=0\\
& &-r_{21}(-v)=r_{12}(v)
\end{eqnarray}
\vspace{1cm}

\noindent{\it\bf Proof:} When $w\longrightarrow 0$, we have following
asympotic properties
\begin{eqnarray*}
& &\sigma (w)=w\sigma'(0)+0(w^3)\ \ \ ,\ \ \
\theta^{(0)}(w)=w\theta^{'(0)}(0)+0(w^3)\\
& &\theta^{(i)}(w)=\theta^{(i)}(0)+w\theta^{'(i)}(0)+0(w^2)\ \ \ ,\ \ 
i\ne 0 \ \ {\rm mod\ \ n}
\end{eqnarray*}
\noindent Then one have ,when $w\longrightarrow 0$
\begin{eqnarray*}
& &\frac{\theta^{'(0)}(0)\sigma(v)}{\sigma'(0)\theta^{(0)}(v)}
\prod_{m=1}^{n-1}\frac{\theta^{(m)}(v)}{\theta^{(m)}(0)}
R^{lk}_{ij}(v)=
w(1-\delta^{l}_{i})\frac{\sigma'(o)}{\sigma(v)}
\prod_{m=1}^{n-1}\frac{\theta^{(m)}(v)}{\theta^{(m)}(0)}
\frac{\theta^{(0)}(v)\theta^{(i-j)}(v)}{\theta^{(l-j)}(v)\theta^{(i-l)}(v)}
+0(w^2)\\
& &\ \ \ \ \
+\delta^{l}_{i}\frac{\sigma'(0)+0(w^2)}{\sigma(v)+w\sigma'(v)+0(w^2)}
\prod_{m=1}^{n-1}\frac{\theta^{(m)}(v)}{\theta^{(m)}(0)}
\frac{\theta^{(0)}(v)(\theta^{(i-j)}(v)+w\theta'^{(i-j)}(v)+0(w^2)}
{\theta^{(i-j)}(v)(\theta'^{(0)}(0)+0(w^2))}\\
& &\ \ \ \ \ \ \  
=\prod_{m=1}^{n-1}\frac{\theta^{(m)}(v)}{\theta^{(m)}(0)}
\frac{\sigma'(0)\theta^{(0)}(v)}{\theta'^{(0)}(0)\sigma(v)}\delta^{l}_{i}
\delta^k_j\\
& &\ \ \ \ \ \
+w\prod_{m=1}^{n-1}\frac{\theta^{(m)}(v)}{\theta^{(m)}(0)}
\{(1-\delta^l_i)\frac{\sigma'(0)\theta^{(0)}(v)\theta^{(i-j)}(v)}
{\sigma(v)\theta^{(l-j)}(v)\theta^{(i-l)}(0)}
+\delta^l_i\delta^k_j\frac{\sigma'(0)\theta^{(0)}(v)}
{\sigma(v)\theta'^{(0)}(0)} (\frac{\theta'^{(i-j)}(v)}{\theta^{(i-j)}(v)}-
\frac{\sigma'(v)}{\sigma(v)})\}\\
& &\ \ \ \ \ \ \ \ \ +0(w^2)
\end{eqnarray*}
\noindent By the definition of classical r-matrix from the quantum one Eq.(15),
we have Eq.(16).The classical Yang-Baxter equation Eq.(17) is the direct
results of the QYBE and the asympotic properties Eq.(18).
The antisymmetric properties of the r-matrix can be derived from the
following relations between the $\theta$-functions
\begin{eqnarray*}
\theta^{(\alpha)}(v)=-e^{2\sqrt{-1}\pi\alpha}\theta^{(-\alpha)}(-v)\ \ \ ,\ \ \
\frac{\theta'^{(\alpha)}(v)}{\theta^{(\alpha)}(v)}=-
\frac{\theta'^{(-\alpha)}(-v)}{\theta^{(-\alpha)}(-v)}\ \ \ {\bf \large\Box}
\end{eqnarray*}
\vspace{1cm}

One can also check that the classical r-matrix $r_{12}(u)$ enjoys in the
$Z_n\otimes Z_n$-symmetric 
\begin{eqnarray*}
r_{12}(v)=(a\otimes a)r_{12}(v)(a\otimes a)^{-1}\ \ \ ,\ \ {\rm for}\ \ a=g,h
\end{eqnarray*}

\section{Review of the elliptic $A_{n-1}$ CM model}
The elliptic $A_{n-1}$ CM model is a system of $n$ one-dimensional
particles interaction by the two-body potential
\begin{eqnarray}
& &V(q_{ij})=\gamma Q(q_{ij})\ \ ,\ \ q_{ij}=q_i-q_j\ \ ,\ \ i,j=1,....,n\\
& &Q(v)-Q(u)=E(u,v)E(u,-v)
\end{eqnarray}
\noindent where $\gamma$ is the coupling constant, $Q(u)$ is a Weierstrass
function and the elliptic function $E(u,v)$ is defined in Eq.(9).In terms of
the canonical variables $\{p_i.q_j\}\ \ (i,j=1,...n)$ with the canonical
Poisson bracket
\begin{eqnarray}
\{p_i,p_j\}=\{q_l,q_k\}=0\ \ ,\ \ \ \{p_i,q_j\}=\delta_{ij}\ \ ,\ \
i,j,k,l=1,...,n
\end{eqnarray}
the Hamiltonian of the system is expressed as
\begin{eqnarray}
H=\sum_{i=1}^{n}p^{2}_i+\sum_{i\ne j}V(q_{ij})
\end{eqnarray}
\noindent The above Hamiltonian with the potential Eq.(19) is known to be
completely integrable[13,22,23,24]. The most effective way to show its
integrability
is to construct the Lax representation for the system. One Lax
pair $(L,M)$ was first found by Krichever[22].The Lax operator (or
L-operator) of Krichever is 
\begin{eqnarray}
L^{i}_{j}(u)=p_i\delta^{i}_{j}+(1-\delta^{i}_{j})\sqrt{\gamma}E(u,q_{ji})
\end{eqnarray}
\noindent where $u$ is spectra parameter and
the motion equation can be rewritten in the Lax form
\begin{eqnarray*}
\frac{d}{dt}L(u)=\{L(u),H\}=[L(u),M(u)]
\end{eqnarray*}
The Hamiltionian defined in Eq.(22) can be rewritten  in terms of the
Poisson-commuting  family $\{trL^{l}(u)\}(l=1,...,n)$,
which forms of enough independent integrals
\begin{eqnarray}
H=tr(L^2(u))+V(u)
\end{eqnarray}
\noindent $V(u)$ does not depend upon the dynamical variables and the
identity Eq.(20) is used.The r-matrix structure of this Lax operator
was given by Sklyanin[30] and Braden et al[10].The fundamental
Poisson bracket of the Lax operator can be described in terms of
r-matrix form[30]
\begin{eqnarray}
\{L_1(u),L_2(v)\}=[r_{12}(u,v),L_1(u)]-[r_{21}(v,u),L_2(v)]
\end{eqnarray}
\noindent and the dynamical r-matrix $r_{12}(u,v)$ is
\begin{eqnarray}
r_{12}(u,v)=a\sum^{n}_{i=1}E^{ii}_{ii}+\sum_{i\ne j}c_{ij}E^{ij}_{ji}
+\sum_{i\ne j}d_{ij}(E^{ii}_{ij}+E^{ji}_{jj})
\end{eqnarray}
\noindent where
\begin{eqnarray}
& &E^{lk}_{ij}=E^{l}_{i}\otimes E^{k}_{j}\ \ ,\ \
a=r^{ii}_{ii}=-\xi(u-v)-\xi(v)\ \ \ \ \ \\
& & c_{ij}=r^{ij}_{ji}=\sqrt{-\gamma}E(u-v,q_{ij})\ \ ,\ \
d_{ij}=r^{ii}_{ij}=r^{ji}_{jj}=\frac{1}{2}\sqrt{-\gamma}E(v,q_{ij})
\end{eqnarray}
\noindent where the elliptic function $\xi(u)$ is defined in Eq.(8).
Sklyanin also shown that the dynamical r-matrix $r_{12}(u,v)$ defined in
Eq.(26) satisfies the dynamical Yang-Baxter equation (or generalized
Yang-Baxter equation)[30]
\begin{eqnarray}
[R^{(123)},L_1] +[R^{(231)},L_2]+[R^{(321)},L_3]=0
\end{eqnarray}
\noindent where
\begin{eqnarray*}
R^{(123)}\equiv r_{(123)}-\{r_{13},L_2\}+\{r_{12},L_3\}
\end{eqnarray*}
\noindent and
\begin{eqnarray*}
r_{(123)}\equiv [r_{12},r_{13}]+[r_{12},r_{23}]-[r_{13},r_{32}]
\end{eqnarray*}
The Jacobi identity of the fundamental Poisson bracket is the results of
Eq.(29).Due to the dynamical properties of the r-matrix $r_{12}(u,v)$ , the
Poisson bracket of L-operator is no longer closed.The quantum version of
Eq.(29) and the generalized (dynamical) Yang-Baxter equation is still not
found except the spin generalization of CM model in which the Gervais-Neveu-
Felder equation was found[4,30].
\section{The ``good" Lax representaion of elliptic $A_{n-1}$ CM model and
its r-matrix}
The L-operator of the elliptic $A_{n-1}$ CM model given by Krichever in
Eq.(23) and corresponding r-matrix $r_{12}(u,v)$ given by Sklyanin in Eq.(26)
leads to some difficulties[30] in the investigation of CM model.This motivate
us to find a ``good" Lax representation of the CM model. As see from
the proposition 1 and its corollary in section 3, this means to find such a
$g(u)$ in Eq.(4) which satisfies Eq.(6).Fortunately, we could find
such a $g(u)$ ,from
which we construct a new L-operator $\stackrel{\sim}{L}(u)$ of the elliptic
$A_{n-1}$ CM model
( This kind L-operator does not always exist for general completely
integrable system).The corresponding r-matrix of $\stackrel{\sim}{L}(u)$
is purely numeric one , and is equal to classical $Z_n$-symmetric r-matrix
In order to compare with the L-operator given by Krichever, we call this
 L-operator found by us as the new Lax operator.

Define
\begin{eqnarray}
& &g(u)=A(u;q)\Lambda(q)\ \ ,\ \ \Lambda(q)^{i}_{j}=h_{i}(q)\delta^{i}_{j}\\
& &h_{j}(q)\equiv h_{j}(q_1,....,q_n)=\frac{1}{\prod_{l\ne i}\sigma(q_{il})}
\nonumber
\end{eqnarray}
\noindent where $A(u;q)^{i}_{j}$ is defined in Eq.(14).Let us construct the
new L-operator $\stackrel{\sim}{L}(u)$
\begin{eqnarray}
& &\stackrel{\sim}{L}(u)=g(u)L(u)g^{-1}(u)\ \ \ \ \ \ \ \ \ \ \\
& &\stackrel{\sim}{M}(u)=g(u)M(u)g^{-1}(u) -(\frac{d}{dt}g(u))g^{-1}(u)
\nonumber
\end{eqnarray}
Then we

\noindent {\bf Proposition 4.} The fundamental Poisson bracket of L-operator
$\stackrel{\sim}{L}(u)$ can be written in the usual Poisson-Lie form with
a purely numerical r-matrix
\begin{eqnarray}
\{\stackrel{\sim}{L}_1(u),\stackrel{\sim}{L}_2(v)\}=[\stackrel{\sim}{r}_{12}
(u-v),\stackrel{\sim}{L}_1(u)+\stackrel{\sim}{L}_2(v)]
\end{eqnarray}
\noindent and the corresponding  r-matrix $\stackrel{\sim}{r}_{12}(u)$ is
a nondynamical one---$Z_n$-symmetric r-matrix defined in Eq.(16).
\vspace{1cm}

\noindent {\it\bf Proof:} The proof is shifted to appendix A.
\vspace{1cm}
The most important properties of this new Lax operator is that the
corresponding r-matrix does not depend upon the dynamical variables.
Consequently, the well-studied theory for the nondynamical system[15]
can be used to study the elliptic $A_{n-1}$ CM model.

The $Z_n$-symmetric r-matrix $\stackrel{\sim}{r}_{12}(u)$ can also be
obtained from the classical dynamical twisting as follows
\begin{eqnarray}
\stackrel{\sim}{r}_{12}(u,v)=g_1(u)g_2(v)r_{12}(u,v)g^{-1}_1(u)g^{-1}_2(v)
+g_2(v)\{g_1(u),L_2(v)\}g^{-1}_1(u)g^{-1}_2(v)
\end{eqnarray}
\noindent up to some matrix which commute with $L_1(u)+L_2(v)$.

The standard Poisson-Lie bracket Eq.(32) of L-operator
$\stackrel{\sim}{L}(u)$ and the numerical r-matrix
$\stackrel{\sim}{r}_{12}(u)$ enjoying in the classical Yang-Baxter equation
Eq(17) and antisymetry Eq.(18) , make it possible to construct the quantum
theory of the elliptic $A_{n-1}$ CM model.Moreover, the numerical r-matrix
$\stackrel{\sim}{r}_{12}(u)$ could provide a mean to construct a seperation
of variables for the elliptic $A_{n-1}$ CM model in the same manner as that
in the case of the integrable magnetic chain[30]. It aslo make it possible
that one can construct the dressing transformation for the model.The dressing
group of this system would be an analogue to the semi-classical limit  of
$Z_n$ Sklyanin algebra.
\section{Discussions}
In this paper, we only construct the nondynamical r-matrix structure for the 
elliptic $A_{n-1}$ Calogero-Moser model. Such a ``good" Lax representation 
for the degenerated case---rational, trigonometric and hyperbolic CM model 
also could be constructed.It is also very interesting to construct such 
a classical dynamical twisting for the Ruijsenaars-Schneider model.

\section*{Appendix}
\subsection*{Appendix A. The proof of Proposition 4}
In this appendix we give the proof the proposition 4 which is main results of
our paper.

Set the classical L-operator $T(u)$ as follow 
\begin{eqnarray*}
T(u)^{i}_{j}=\sum_{k}A(u;q)^{i}_{k}A^{-1}(u;q)^{k}_{j}p_{k}
-s(\partial_{u}A(u;q))^{i}_{k}A^{-1}(u;q)^{k}_{j}
\end{eqnarray*}
\noindent where $\{p_{k}\}$ is the classical moment which is conjugated
with $\{q_k\}$ and $\{p_i,q_j\}$ satisfies the canonical Poisson bracket
Eq.(21).

{\bf Lemma 1.} The classical operator $T(u)$  has the standard Poisson-Lie
bracket
\begin{eqnarray}
\{T_1(u),T_2(v)\}=[\stackrel{\sim}{r}_{12}(u-v),T_1(u)+T_2(v)]
\end{eqnarray}
\noindent where r-matrix $\stackrel{\sim}{r}_{12}(u)$ is the $Z_n$-symmetric
r-matrix defined in Eq.(16).
\vspace{1cm}

\noindent {\it\bf Proof :} When $w\longrightarrow 0$ , the difference quantum
L-operator $\widehat{T}(u)$ has the following asympotic properties
\begin{eqnarray*}
& &\widehat{T}(u)^{i}_{j}=\sum_{k}A(u;q)^{i}_{k}A^{-1}(u;q)^{k}_{j}-
w\sum_{k}A(u;q)^{i}_{k}A^{-1}(u;q)^{k}_{j}\frac{\partial}{\partial q_{k}}\\
& &\ \ \ \ \ \ \
+sw\sum_{k}(\partial_{u}A(u;q))^{i}_{k}A^{-1}(u;q)^{k}_{j}+ 0(w^2)\\
& &\ \ \ \equiv \delta^{i}_{j}-w\widehat{T}^{(1)}(u)^{i}_{j}+0(w^2)
\end{eqnarray*}
\noindent where
\begin{eqnarray*}
\widehat{T}(u)^{i}_{j}=\sum_{k}A(u;q)^{i}_{k}A^{-1}(u;q)^{k}_{j}
\frac{\partial}{\partial q_{k}}-s\sum_{k}(\partial_{u}A(u;q))
^{i}_{k}A^{-1}(u;q)^{k}_{j}
\end{eqnarray*}
\noindent From the QYBE Eq.(13), we have 
\begin{eqnarray*}
[\widehat{T}^{(1)}_1(u),\widehat{T}^{(1)}_2(v)]=
[\stackrel{\sim}{r}_{12}(u-v),\widehat{T}^{(1)}_1(u)+\widehat{T}^{(1)}_2(v)]
\end{eqnarray*}
\noindent If we use $p_k$ instead of the differential
$\frac{\partial}{\partial q_{k}}$ and the classical L-operator $T(u)$
instead of $\widehat{T}^{(1)}(u)$ , we have the Eq.(34).
\hspace{3cm}${\bf \large\Box}$

\noindent {\bf Lemma 2.} The $T(u)$ can be explictly written as follows
\begin{eqnarray}
T(u)^{i}_{j}=\sum_{k,k'}\stackrel{\sim}{g}(u)^{i}_{k}\overline{T}(u)^{k}_{k'}
\stackrel{\sim}{g}^{-1}(u)^{k'}_j
\equiv \{(A(u;q)\Lambda(q))\overline{T}(u)(A(u;q)\Lambda(q))^{-1}\}^{i}_{j}
\end{eqnarray}
\noindent where $\overline{T}(u)$ and $\Lambda(q)$ are
\begin{eqnarray}
& &\overline{T}(u)^{i}_{j}=(p_{i}-\frac{\partial}{\partial q_k}
ln\Delta^{\frac{s}{n}}(q))\delta^{i}_{j}+\sqrt{\gamma}(1-\delta^{i}_{j})
E(u;q_{ji})
\nonumber\\
& & \Delta(q)=\prod_{i<j}\sigma(q_{ij})\ \ \ ,\ \ \
{\rm coupling \ \ constant\ \ } \gamma=(-\frac{s\sigma'(0)}{n})^{2}
\ \ \ ,\ \ \
\Lambda(q)^{i}_{j}=\frac{1}{\prod_{l\ne i}\sigma(q_{il})}\delta^{i}_{j}
\end{eqnarray}
\vspace{1cm}

\noindent {\it\bf Proof:}
In order to calculate the matrix element $(\partial_{u}(A(u;q))A^{-1}(u;q)$ ,
we first consider
\begin{eqnarray*}
A(u+w;q)A^{-1}(u;q)=A(u;q)[A^{-1}(u;q)A(u+w;q)]A^{-1}(u;q)
\end{eqnarray*}
\noindent From the defintion of $A(u;q)^{i}_{j}$ Eq.(14) and the determinat
formula of Vandermonde type
\begin{eqnarray*}
det[\theta^{(j)}(u_k)]=const.\times \sigma(\frac{1}{n}\sum_{k}u_{k}
-\frac{n-1}{2})
\prod_{1\leq j<k\leq n}\sigma(\frac{u_k-u_j}{n})
\end{eqnarray*}
\noindent where the $const.$ does not depend upon $\{u_{k}\}$, we have
\begin{eqnarray*}
[A^{-1}(u;q)A(u+w;q)]^{i}_{j}=\sum_{k}A^{-1}(u;q)^{i}_{k}A(u+w;q)^{k}_{j}
=\frac{\sigma(\frac{w}{n}+u+q_{ji})}{\sigma(u)}\prod_{k\ne i}
\frac{\sigma(\frac{w}{n}+q_{jk})}
{\sigma(q_{ik})}
\end{eqnarray*}
\noindent Then, we have
\begin{eqnarray*}
& &(A^{-1}(u;q)\partial_{u}A(u;q))^{i}_{j}=\frac{\partial}{\partial w}
\left\{\frac{\sigma(\frac{w}{n}+u+q_{ji})}{\sigma(u)}\prod_{k\ne i}
\frac{\sigma(\frac{w}{n}+q_{jk})}
{\sigma(q_{ik})}\right\}|_{w=0}\\
& &\ \ \ =\frac{1}{n}
\{\frac{\sigma'(u)}{\sigma(u)}\delta^{i}_{j}+\frac{\sigma(u+q_{ji})}
{\sigma(u)}(\delta^{i}_{j}\sum_{k\ne i}\frac{\sigma'(q_{ik})}
{\sigma(q_{ik})}+(1-\delta^{i}_{j})\sigma'(0)\frac{\prod_{k\ne i,j}
\sigma(q_{jk})}
{\prod_{k\ne i}\sigma(q_{ik})})\}\\
& &\ \ \ =\frac{1}{n}
\{(\frac{\sigma'(u)}{\sigma(u)}+\sum_{k\ne i}\frac{\sigma'(q_{ik})}
{\sigma(q_{ik})})\delta^{i}_{j}+(1-\delta^{i}_{j})
\frac{\sigma'(0)\sigma(u+q_{ji})}{\sigma(u)\sigma(q_{ji})}
\frac{\prod_{k\ne j}\sigma(q_{jk})}{\prod_{k\ne i}\sigma(q_{ik})}\}\\
& &\ \ \ =\frac{1}{n}
\{(\frac{\sigma'(u)}{\sigma(u)}+\frac{\partial}{\partial q_{j}}
(ln\Delta(q)))\delta^{i}_{j}-(1-\delta^{i}_{j})(-\sigma'(0))E(u;q_{ji})
\frac{\prod_{k\ne j}\sigma(q_{jk})}{\prod_{k\ne i}\sigma(q_{ik})}\}\\
& &\ \ \ =
\frac{1}{\prod_{k\ne i}\sigma(q_{ik})}
\{(\frac{\sigma'(u)}{n\sigma(u)}+\frac{\partial}{\partial q_{j}}
(ln\Delta^{\frac{1}{n}}(q)))\delta^{i}_{j}-(1-\delta^{i}_{j})
(-\frac{\sigma'(0)}{n})E(u;q_{ji})\}
\prod_{k\ne i}\sigma(q_{ik})
\end{eqnarray*}
\noindent Subsituting $A^{-1}(u;q)\partial_{u}A(u;q)$ into the defintion of
$T(u)$ , we have
\begin{eqnarray*}
& &T(u)^{i}_{j}=\sum_{k}A(u;q)^{i}_{k}p_{k}A^{-1}(u;q)^{k}_{j}+
-s\sum_{k,k'}A(u;q)^{i}_{k}(A^{-1}(u;q)\partial_{u}A(u;q))^{k}_{k'}
A^{-1}(u;q)^{k'}_{j}\\
& &\ \ \ =
(A(u;q)\Lambda(q))^{i}_{k}
\{(p_{k}-\frac{s\sigma'(u)}{n\sigma(u)}+\frac{\partial}{\partial q_{k}}
(ln\Delta^{\frac{s}{n}}(q)))\delta^{k}_{k'}+(1-\delta^{k}_{k'})
\sqrt{\gamma}E(u;q_{ji})\}
(\Lambda(q)^{-1})A(u;q)^{-1})^{k'}_{j}
\end{eqnarray*}
\hspace{4cm} $\ \ \ \ \ \ \ \ \ \ \ \ \ \ \ \ \ {\bf\large \Box}$

\noindent We consider a map:
\begin{eqnarray}
\left\{\begin{array}{l}
p_{i}\longrightarrow p_{i}-\frac{\partial}{\partial q_i}(ln\Delta
^{\frac{s}{n}}(q))\\
q_{i}\longrightarrow q_{i}
\end{array}
\right.
\end{eqnarray}
\noindent{\bf Lemma 3.} The map defined in Eq.(37) is a Poisson map[1] (
or a canonical transformation).
\vspace{1cm}

\noindent{\it\bf Proof:}The Lemma 3 can be proven from considering
the symplectic two-form
\begin{eqnarray*}
& &\sum_{i}d(p_{i}-\frac{\partial}{\partial q_{i}}
ln\Delta^{\frac{s}{n}}(q))\wedge dq_{i}
= \sum_{i}dp_{i}\wedge dq_{i}+\sum_{ij}(\frac{\partial^2}
{\partial q_i\partial q_j}ln\Delta^{\frac{s}{n}}(q))dq_{i}\wedge dq_{j}\\
& &\ \ \ \ \ =\sum_{i}dp_{i}\wedge dq_{i}\ \ \ \ {\bf \large\Box}
\end{eqnarray*}
\noindent Since the Poisson bracket is invariant under the Poisson map[1],
 we could have proposition 4 from Lemma 1 and Lemma 3.

\newpage
\section*{References}
\begin{enumerate}

\item Arnold,V.I. : Mathematical methods of classical mechanics , Springer
Verlag (1978).

\item Arutyunov,G.E., Frolov,S.A., Medredev,P.B. : Elliptic
Ruijsenaars-Schneider model via the Poisson reducation of the affine
Heiserberg Double, {\it hep-th/9607170}(1996); Elliptic Ruijsenaars-Schneider
model form the cotangent bundle over the two-dimensional current group,
{\it hep-th/9608013}(1996).

\item Avan,J. : Classical dynamical r-matrices for Calogero-Moser systems
and their generalizations,

\item Avan,J., Babelon,O., Billey,E.: {\it The Gervais-Neveu-Felder 
equation and the quantum Calogero-Moser systems}, hep-th/9505091 (1995)
{\it Comm. Math. Phys.}{ \bf 178}, 281(1996).

\item Avan,J., Talon,T.: {\it Phys. Lett.} {\bf B 303}, 33(1993).

\item Babelon,O., Bernard,D.:{\it Phys. Lett.} {\bf B 317},363 (1993).

\item Babelon,O., Viallet,C.M. : {\it Phys. Lett.} {\bf B237}, 411(1989).

\item Belavin,A.A. : {\it Nucl. Phys.} {\bf B180}, 189(1981).

\item Belavin,A.A., Drinfeld,V.G. : Triangle equation and simple Lie
algebras, {\it Soviet Sci. reviews}, {\bf Sect.C 4}, 93(1984).

\item Braden,H.W., Suzaki,T.: {\it Lett. Math. Phys.} {\bf Vol.30}, 147(1994).

\item Braden,H.W, Andrew,N.W.Hone : Affine Toda solitons and systems
of Calogero-Moser type, {\it hep-th/9603178}(1996).

\item Braden,H.W., Corrigan,E., Dorey,P.E., Sasaki,R. : {\it Nucl. Phys.}
{\bf B338}, 689(1990); {\it Nucl. Phys.} {\bf B356}, 469(1991).

\item Calogero,F. : {\it Lett. Nuovo. Cim.} {\bf 13}, 411 (1975);
{\it Lett.Nuovo. Cim.} {\bf 16}, 77 (1976).

\item Etingof,P., Varchenko,A. :Geometry and classification of solutions 
of the dynamical Yang-Baxter equation, {\it q-alg/9703040}.

\item Faddeev,L.D., Takhtajan,L. :Hamiltonian methods in the theory of
solitons , Springer Verlag (1987).

\item Gorsky,A., Nekrasov,N.: {\it Nucl. Phys.} {\bf B414}, 213(1994);
{\it Nucl. Phys.} {\bf B436}, 582(1995).

\item Hasegawa,K., :Ruijsenaars' commuting difference operators as commuting
transfer matrices, {\it q-alg/9512029}; {\it Jour. Math. Phys. }{\bf 35},
6158 (1994).

\item Hou,B.Y., Yang,W.L., :The nondynamical r-matrix structure of the 
elliptic Calogero-Moser model (n=2), {\bf Preprint :} {\it IMPNWU-960810}.

\item Hou,B.Y., Shi,K.J., Yang,Z.X.: J. Phys. {\bf A 26} 4951(1993).

\item Hou,B.Y., Wei,H. : {\it J. Math. Phys.} {\bf 30}, 2750(1989).

\item Jimbo,M., Miwa,T., Okado,M. :{\it Nucl. Phys.} {\bf B300}, 74(1988).

\item Krichever,I.M.: {\it Func. Annal. Appl.} {\bf 14},  282(1980). 

\item Moser,J.: {\it Adv. Math.} {\bf 16}, 1(1975).

\item Olshanetsky,M.A., Perelomov,A.M.,  {\it Phys. Rep.} {\bf 71}, 313(1981).

\item Quano,Y.H., Fujii,A., :{\it Mod. Phys. Lett.} {\bf A6}, 3635 (1991).

\item Richey,M.P., Tracy,C.A. : {\it J. Stat. Phys. }{bf 42}, 311(1986);
Tracy,C.A.: {\it Physica }{\bf D16}, 203(1985).

\item Ruijsenaars,S.N.M., Schneider,H.: {\it Ann. Phys. }
{\bf Vol.170}, 370(1986).  

\item Ruijsenaars,S.N.M.: {\it Comm. Math. Phys.} {\bf Vol.115}, 127(1988). 

\item Schiffman,O. : On classification of dynamical r-matrices,
{\it q-alg/9706017}(1997).

\item Sklyanin,E.K.:{\it Funct. Anal. Appl.} {\bf 16}, 263 (1982);
{\it Funct.Anal. Appl.} {\bf 17}, 320(1983);{\it Comm.Math.Phys.}
{\bf Vol.150}, 181(1992); {\it Dynamical r-matrices for
the elliptic Calogero-Moser model}, hep-th/9308060 (1993).

\item Suris,Yuri B.: {\it Why are the rational and hyperbolic 
Ruijsenaars-Schneider hierarchies governed by the the same R-operators as 
the Calogero-Moser ones}, hep-th/9602160 (1996).

\end{enumerate}
\end{document}